# The Cohomology of the Lyons Group and Double Covers of Alternating Groups

Alejandro Adem[†]　　Dikran Karagueuzian[‡]
Department of Mathematics
University of Wisconsin, Madison WI 53706

R. James Milgram[†]
Department of Mathematics
Stanford University, Stanford CA 94305

Kristin Umland
Department of Mathematics
University of New Mexico, Albuquerque NM 87131

## §0 INTRODUCTION

From the Classification Theorem for Finite Simple Groups we know that there are precisely eight sporadic simple groups which have $(\mathbb{Z}/2)^4$ as their largest elementary abelian 2–subgroup (i.e. they have rank equal to four at the prime $p = 2$): the two Mathieu groups $M_{22}$, $M_{23}$; the two Janko groups $J_2$, $J_3$; the McLaughlin group $McL$; the Lyons group $Ly$; and two groups with much larger 2-Sylow subgroups $HS$ (Higman–Sims) and $Co_3$ (one of Conway's groups). In previous work the mod 2 cohomology rings of $M_{22}$ [AM2], $M_{23}$ [M], $J_2$ and $J_3$ [CMM], and $McL$ [AM3] were computed. One of the main objectives of this paper is to determine the cohomology ring $H^*(Ly; \mathbb{F}_2)$.

Even though $Ly$ is a large and complicated group (it is the largest of those in the list of above, having order equal to $51,765,179,004 \times 10^6$), the answer turns out to be surprisingly simple. We show that $Ly$ has two conjugacy classes of maximal 2-elementary subgroups, one $2^4$ with Weyl group $A_7$, and one $2^3$ with Weyl group $L_3(2)$, and we prove that these subgroups detect $H^*(Ly)$. In [AM4] we determined the ring of invariants

$$H^*(2^4)^{A_7} = \mathbb{F}_2[d_8, d_{12}, d_{14}, d_{15}](1, a_{18}, a_{20}, a_{21}, a_{24}, a_{25}, a_{27}, a_{45})$$

while it is well–known that the invariants $H^*(2^3)^{L_3(2)}$ are $\mathbb{F}_2[d_4, d_6, d_7]$, the Dickson algebra. From this we obtain

[†] Partially supported by a grant from the NSF.
[‡] Supported by an NSF Postdoctoral Fellowship.

**Theorem 7.3**

There is a short exact sequence of rings

$$0 \longrightarrow H^*(Ly; \mathbb{F}_2) \longrightarrow H^*(2^4)^{A_7} \oplus \mathbb{F}_2[d_4^2, d_6^2, d_7](1, d_4d_7, d_6d_7, d_4d_6d_7) \longrightarrow \mathbb{F}_2[h_8, h_{12}] \longrightarrow 0$$

with $h_8$ the image of $d_8$ in the $A_7$ invariants and $d_4^2$ in the $L_3(2)$ invariants, while $h_{12}$ is the image of $d_{12}$ and $d_6^2$ respectively. ∎

Thus we can write

$$H^*(Ly, \mathbb{F}_2) \cong H^*(2^4)^{A_7} \oplus d_7 \mathbb{F}_2[d_4, d_6, d_7]$$

where the role of $d_7$ and $d_4 d_7$ is to have $\mathbb{Z}/8$–Bocksteins to $d_8$ and $d_{12}$ respectively. In the course of obtaining this result we determine the cohomology of a number of groups which arise in other contexts and perhaps our results for these further groups are even more intriguing.

The double covers of the alternating groups $2A_n$ arise in studying the two connected cover of the infinite loop space $Q(S^0)$, which can actually be identified with the plus construction $B_{2A_\infty}^+$. The mod 2 cohomology of $B_{2A_\infty}$ is easily determined since $H^*(B_{A_\infty}, \mathbb{F}_2)$ is known as a polynomial algebra on a specified but infinite set of generators, and the elements $\sigma_2$, $Sq^1(\sigma_2)$, $Sq^2Sq^1(\sigma_2)$, ..., $Sq^{2^n}Sq^{2^{n-1}}\cdots Sq^2Sq^1(\sigma_2)$, ..., form a regular sequence where $\sigma_2$ is the non-zero element in dimension two. However, the cohomology of the individual $2A_n$'s is much more complex and the inclusions $2A_n \hookrightarrow 2A_{n+1}$ tend to often have very large cokernels in cohomology.

By studying these cokernels it appears that they give information related to the structure of periodicity operators in homotopy theory, and (as one sees in the proof of 4.2 where the exact structure of the cokernel for $2A_8$ is seen to be equivalent to the determination of the ring of invariants) also appear to play an intriguing role in modular invariant theory.

Moreover, $2A_{11}$ is a maximal, odd index subgroup of $Ly$ while the Sylow 2-subgroup of $2A_8$ is the Sylow 2-subgroup of $M_{22}$, $M_{23}$ and $McL$. Our route to the cohomology of $Ly$ determines the cohomology of the following groups in succession

$$2A_8 \hookrightarrow 2S_8 \hookrightarrow 2A_{10} \hookrightarrow 2A_{11}$$

where the last inclusion is a 2-local equivalence. We now describe how this goes.

It is well–known that $H^2(A_n, \mathbb{F}_2) \cong \mathbb{F}_2$, hence for each of them there exists a unique non–trivial double cover

$$1 \to \mathbb{Z}/2 \to 2A_n \to A_n \to 1.$$

An interesting group–theoretic fact is that for $n \leq 7$, $2A_n$ does not contain $\mathbb{Z}/2 \times \mathbb{Z}/2$ as a subgroup. Hence $H^*(2A_n, \mathbb{F}_2)$ is periodic in this range, and in fact they are all isomorphic to $\mathbb{F}_2[u_4] \otimes \Lambda(x_3)$, a polynomial algebra on a 4–dimensional generator tensored with an exterior algebra on a 3–dimensional class. This pleasant behaviour breaks down when $n = 8$ and in fact the calculation of $H^*(2A_8, \mathbb{F}_2)$ has been an open problem for some time. Attempts to use the description of $H^*(A_8, \mathbb{F}_2)$ provided in [AMM] have been unsuccessful due to the *a priori* complicated nature of the spectral sequence associated to the extension above. We show:

**Theorem 3.3**

The mod 2 cohomology of $2A_8$ can be described by a long exact sequence

$$0 \to \mathbb{F}_2[d_4, d_8](u_3, u_7, u_9) \to H^*(2A_8, \mathbb{F}_2) \to H^*(E)^{L_3(2)} \oplus H^*(F)^{L_3(2)} \to \mathbb{F}_2[d_4, d_8] \to 0$$

where $E$ and $F$ are rank 4 elementary abelian subgroups in $2A_8$, and

$$H^*(E)^{L_3(2)} = H^*(F)^{L_3(2)} = \mathbb{F}_2[d_4, d_6, d_7, d_8](1, a_8, a_9, a_{10}, a_{11}, a_{12}, a_{13}, a_{21}). \blacksquare$$

In the expression above, $E$ and $F$ represent the two conjugacy classes of maximal elementary abelian subgroups. The term on the left is the *radical* (i.e. the nilpotent elements) in $H^*(2A_8)$. The term on the right represents the double image classes. We use the above to compute $H^*(2S_8, \mathbb{F}_2)$, where $2S_8$ is the double cover of $S_8$ which extends the double cover of the alternating group. Unexpectedly, it turns out that the elementary abelian subgroups *detect* the cohomology of this double cover. More precisely we have

**Theorem 5.4**

The cohomology of $2S_8$ is determined by the following exact detection sequence

$$0 \to H^*(2S_8, \mathbb{F}_2) \to \mathbb{F}_2[d_4, d_6, d_7, d_8](1, a_8, a_9, a_{10}, a_{11}, a_{12}, a_{13}, a_{21})$$
$$\oplus \mathbb{F}_2[w_1, d_2^2, d_4^2](1, d_3, d_3 d_4, t(t+w)d_3 d_4) \to \mathbb{F}_2[h_4, h_8] \to 0$$

where $h_4, h_8$ correspond to the double image classes $(d_4, d_2^2)$ and $(d_8, d_4^2)$ respectively. $\blacksquare$

One can deduce from this that the spectral sequence for the extension $2A_8 \triangleleft 2S_8$ collapses at $E_2$; the notation reflects this fact. Here we have once again two conjugacy classes of maximal elementary abelian subgroups, one of rank 4 (obtained by fusing the 2 representatives for $2A_8$) and one of rank 3. As before, the term on the right represents the double image classes.

Noting that $2S_8 \subset 2A_{10}$ is an odd index subgroup, we infer that $H^*(2A_{10}, \mathbb{F}_2)$ is also detected on elementary abelian subgroups. The maximal elementary abelian subgroups are weakly closed in $2S_8 \subset 2A_{10}$ and hence, using the Cárdenas–Kuhn Theorem (see [AM1]) the image of the restriction maps can be obtained as the intersection of the restriction from $H^*(2S_8, \mathbb{F}_2)$ with the invariants under the appropriate Weyl groups: still $L_3(2)$ for the $2^4$ and now $S_4$ (instead of $D_8$) for the $2^3$. In fact it can be described as follows.

**Theorem 6.1**

The cohomology of $2A_{10}$ is determined by the following exact detection sequence

$$0 \to H^*(2A_{10}, \mathbb{F}_2) \to H^*(2^4)^{L_3(2)} \oplus \mathbb{F}_2[d_3, d_2^2, d_4^2](1, d_2 d_3, d_3 d_4, d_2 d_3 d_4) \to \mathbb{F}_2[h_4, h_8] \to 0$$

where $h_4, h_8$ correspond to the double image classes $(d_4, d_2^2)$ and $(d_8, d_4^2)$ respectively. $\blacksquare$

To obtain $H^*(Ly, \mathbb{F}_2)$ from the expression above, we simply need to take $A_7$ invariants and compute the intersection of the image of $H^*(2A_{10})$ in $H^*(2^3)$ with the Dickson Algebra of $L_3(2)$–invariants.

It would seem that Theorem 5.4 should admit an extension to double covers of $S_{2^i}$ for $i \geq 4$. These groups carry important homotopy–theoretic data in their cohomology and hence deserve further attention. We will defer this to a sequel.

This paper is organized as follows: in §1 we discuss preliminaries; in §2 we determine the double image classes for the cohomology of $2A_8$; in §3 we determine the radical of $H^*(2A_8)$; in §4 we consider the spectral sequence associated to the central extension $1 \to \mathbb{Z}/2 \to 2A_8 \to A_8 \to 1$; in §5 we describe facts about the 2–Sylow subgroup of $Ly$ and determine the cohomology of $2S_8$; in §6 we determine the cohomology of $2A_{10}$; and in §7 we calculate $H^*(Ly, \mathbb{F}_2)$.

We are grateful to L. Evens, M. Isaacs and S. Priddy for helpful comments. We thank Gregor Kemper for the use of a beta test version of his invariant theory package, and for the method of argument with which we determined $H^*(2^4)^{L_3(2)}$. Some of the calculations in this paper were first performed or subsequently checked using Magma. We thank the Magma project for their assistance in the use of this valuable tool. We thank R. Lyons for providing a purely group-theoretic proof of Corollary 5.5. Also we should point out that the first indication that such detection results might be possible occurs in [U] where, besides determining the Weyl group of the new $2^3$ in $Ly$, it was shown that the element $u_7$ when looked at in $H^*(Ly, \mathbb{F}_2)$ was not nilpotent.

Throughout this paper, coefficients will be assumed taken in the field $\mathbb{F}_2$ with two elements, so they are suppressed.

## §1. PRELIMINARIES

We begin by discussing the mod 2 cohomology of the double covers $2A_n$ for $n \leq 7$. In this range all these groups have periodic cohomology rings at the prime 2, and they can be easily computed. In fact we have the isomorphisms $2A_4 \cong SL_2(\mathbb{F}_3)$, $2A_5 \cong SL_2(\mathbb{F}_5)$, $2A_6 \cong SL_2(\mathbb{F}_9)$, while a simple double coset argument shows that $2A_7$ is mod 2 cohomologous to $2A_6$. Note that for the same reason $2A_4$ is mod 2 cohomologous to $2A_5$. We state their cohomology, leaving details to the reader (or see [AM1], Ch.IV):

**Proposition 1.0**
  If $G = 2A_n$, where $n \leq 7$, then $H^*(G) \cong \mathbb{F}_2[v_4] \otimes E(y_3)$. ∎

This remarkably simple behaviour breaks down rather dramatically for $n = 8$. Even though the mod 2 cohomology of $A_8$ was computed in [AMM] some time ago, the computation for its double cover poses serious technical difficulties. Our strategy will be to determine $H^*(2A_8)$ using techniques developed to study the cohomology of certain sporadic simple groups. This calculation in turn will allow us to compute the cohomology of $2A_{10}$ and from there that of the sporadic simple group $Ly$.

To begin we provide background information needed to compute the mod 2 cohomology of the group $G = 2A_8$. First we recall that $G$ has precisely two conjugacy classes of maximal elementary abelian 2–subgroups, both of rank 4, which from now on we will denote by $E$ and $F$ (see [AM3]). These are constructed as follows. First consider the regular representation of $(\mathbb{Z}/2)^3 = V_3 \subset S_8$. This representation has Weyl group $Aut((\mathbb{Z}/2)^3) = GL_3(2)$. (As is always the case for the regular representation $G \hookrightarrow S_{|G|}$, we have $N(G)/G = Out(G)$.) However, restricting $V_3$ to $A_8$ gives two non-conjugate subgroups $V_3$, $V_3'$ in $A_8$ which fuse in $S_8$, and both $V$, $V'$ lift to $(\mathbb{Z}/2)^4$'s in $2A_8$. In fact, if $H = Syl_2(2A_8)$, then it contains precisely two copies of $(\mathbb{Z}/2)^4$. The subgroups $E$ and $F$ intersect in a rank 2 subgroup $V = E \cap F$ and

**Lemma 1.1**:

The normalizers of $E$, $F$ in $2A_8$ are split extensions of the form $E : L_3(2)$ and $F : L_3(2)$, which are isomorphic since the two $(\mathbb{Z}/2)^4$'s fuse in $2S_8 = (2A_8):2$. In particular both of these maximal 2-tori have $L_3(2)$ as their Weyl group. ∎

The action of these matrix groups on $E$, $F$ can be written down explicitly as follows: $L_3(2)$ is the subgroup of $L_4(2)$ generated by the following two matrices

$$A = \begin{pmatrix} 1 & 1 & 1 & 0 \\ 1 & 0 & 1 & 0 \\ 1 & 1 & 0 & 0 \\ 1 & 1 & 0 & 1 \end{pmatrix} \qquad B = \begin{pmatrix} 0 & 0 & 1 & 0 \\ 0 & 1 & 0 & 0 \\ 1 & 1 & 1 & 0 \\ 1 & 1 & 0 & 1 \end{pmatrix}$$

with simple components the usual action of $L_3(2)$ on $(\mathbb{Z}/2)^3$ and the trivial action on $\mathbb{Z}/2$. But the representation is *non-split* and this will be reflected in the structure of the ring of invariants $H^*((\mathbb{Z}/2)^4)^{L_3(2)}$.

Recall that the $\mod(2)$ cohomology of $2^4$, $H^*(2^4)$, is the polynomial ring on four one dimensional variables, $\mathbb{F}_2[x_1, y_1, z_1, w_1]$ and the action of $L_3(2)$ is given by the *transposes* of $A$ and $B$ above. Consequently, the subring $\mathbb{F}_2[x_1, y_1, z_1]$ is invariant under the action and the Dickson algebra $\mathbb{F}_2[x_1, y_1, z_1]^{L_3(2)} = \mathbb{F}_2[d_4, d_6, d_7] \subset \mathbb{F}_2[x_1, y_1, z_1, w_1]^{L_3(2)}$. In fact, by a standard computer assisted computation of invariants, we have

**Proposition 1.2**:

The ring of invariants $\mathbb{F}_2[x_1, y_1, z_1, w_1]^{L_3(2)}$ can be described as

$$\mathbb{F}_2[d_4, d_6, d_7, d_8](1, a_8, a_9, a_{10}, a_{11}, a_{12}, a_{13}, a_{21})$$

where the subscript denotes the dimensions of the generators. The ring is freely generated over the polynomial subring spanned by $d_4, d_6, d_7, d_8$ and furthermore we have $a_9 = Sq^1 a_8$, $a_{11} = Sq^2 Sq^1 a_8$, $a_{12} = Sq^2 a_{10}$, $a_{13} = Sq^1 Sq^2 a_{10}$ and $a_{21} = a_{10} a_{11}$.

**Proof**: Write $W \cong L_3(2)$. It was shown in [AM2] that for a 2-Sylow subgroup $P \subset W$, the invariant ring $H^*(E)^P$ is Cohen-Macaulay; it therefore follows (see [B]) that $H^*(E)^W$ is Cohen-Macaulay also.

Using the modular invariant theory package developed by Gregor Kemper [K] in Maple, one obtains a system of primary invariants for $H^*(E)^W$ having degrees $4, 6, 7,$ and $8$. Suitable choices for these invariants are the Dickson invariants $d_4, d_6, d_7$, described above, and $d_8$, the first Dickson invariant for $H^*((\mathbb{Z}/2)^4)$,

$$d_8 = w_1^8 + w_1^4 d_4 + w_1^2 d_6 + w_1 d_7 + d_4^2.$$

The "folk-theorem" given by Kemper in [K] now tells us how many secondary invariants are necessary to generate the invariant ring. The number of secondary invariants is

$$\frac{(4 \cdot 6 \cdot 7 \cdot 8)}{|W|} = 8.$$

To obtain the secondary invariants one determines, using Kemper's package, the dimensions of the spaces of homogeneous invariants in degrees $\leq 13$, which are given in the following table:

| degree : | 0 | 1 | 2 | 3 | 4 | 5 | 6 | 7 | 8 | 9 | 10 | 11 | 12 | 13 |
|---|---|---|---|---|---|---|---|---|---|---|---|---|---|---|
| dimension : | 1 | 0 | 0 | 0 | 1 | 0 | 1 | 1 | 3 | 1 | 2 | 2 | 5 | 3 |

Through degree 13, the polynomial algebra generated by the primary invariants has the following dimensions:

| degree : | 0 | 1 | 2 | 3 | 4 | 5 | 6 | 7 | 8 | 9 | 10 | 11 | 12 | 13 |
|---|---|---|---|---|---|---|---|---|---|---|---|---|---|---|
| dimension : | 1 | 0 | 0 | 0 | 1 | 0 | 1 | 1 | 2 | 0 | 1 | 1 | 3 | 1 |

Thus, the first degree in which there is a shortfall is 8. Therefore we must have a secondary invariant $a_8$ of degree 8. Writing $R$ for $\mathbb{F}_2[d_4, d_6, d_7, d_8]$, we can compute the graded dimension in degrees 0 through 13 of the $R$-module $R + R \cdot a_8$, which is of course an $R$-submodule of the invariant ring:

| degree : | 0 | 1 | 2 | 3 | 4 | 5 | 6 | 7 | 8 | 9 | 10 | 11 | 12 | 13 |
|---|---|---|---|---|---|---|---|---|---|---|---|---|---|---|
| dimension : | 1 | 0 | 0 | 0 | 1 | 0 | 1 | 1 | 3 | 0 | 1 | 1 | 4 | 1 |

We see now that the next shortfall is in dimension 9; therefore there must be a secondary invariant of degree 9, call it $a_9$. Again, we can compute the graded dimension of $R + R \cdot a_8 + R \cdot a_9$, we obtain

| degree : | 0 | 1 | 2 | 3 | 4 | 5 | 6 | 7 | 8 | 9 | 10 | 11 | 12 | 13 |
|---|---|---|---|---|---|---|---|---|---|---|---|---|---|---|
| dimension : | 1 | 0 | 0 | 0 | 1 | 0 | 1 | 1 | 3 | 1 | 1 | 1 | 4 | 2 |

A comparison of this table with the table of dimensions of the full invariant ring, in just the manner of the previous two examples, shows that there must be secondary invariants in degrees $10, 11, 12$, and $13$, which we will call $a_{10}, a_{11}, a_{12}$, and $a_{13}$.

Finally we may check that $W$ contains no "pseudo-reflection", using Maple, and that we may therefore apply the Carlisle-Kropholler conjecture [B] to obtain the fact that the final secondary invariant must have degree 21. We call this final secondary invariant $a_{21}$ and note that the computer tells us that we may in fact take $a_{21} = a_8 \cdot a_{13}$. Furthermore, using the package's ability to compute Steenrod operations, we note that we may take $a_9 = Sq^1 a_8$, $a_{11} = Sq^2 a_9$, $a_{12} = Sq^2 a_{10}$, and $a_{13} = Sq^1 a_{12}$.

Finally we note that this implies the Poincare series of the invariant ring is

$$\frac{1 + t^8 + t^9 + t^{10} + t^{11} + t^{12} + t^{13} + t^{21}}{(1 - t^4)(1 - t^6)(1 - t^7)(1 - t^8)}$$

and that the invariants are of the form

$$\mathbb{F}_2[d_4, d_6, d_7, d_8](1, a_8, a_9, a_{10}, a_{11}, a_{12}, a_{13}, a_8 a_{13}). \quad \blacksquare$$

If we restrict ourselves to the 2–Sylow subgroup $H$, then the corresponding Weyl group for $E$ and $F$ is the dihedral group $D_8$. In [AM2] it was shown that the cohomology $H^*(H, \mathbb{F}_2)$ maps onto the $D_8$ invariants in $H^*(E)$ and $H^*(F)$. Clearly the image of the restriction from $H^*(2A_8)$ will map into the $L_3(2)$ invariants, however we in fact claim that it maps onto them in both cases. To show this we first recall a basic notion useful in cohomology calculations. A triple $(Q, H, K)$ of groups with $K \subset H \subset Q$ are said to form a *weakly closed system* if any subgroup in $Q$ conjugate to $K$ is already conjugate to $K$ in $H$. The Cárdenas-Kuhn Theorem (see [AM1]) asserts that if $(Q, K, V)$ is a weakly closed system with $V = (\mathbb{Z}/2)^r$ and $[Q : K]$ odd, then the image of the restriction map from $H^*(Q)$ to $H^*(V)$ is precisely the image of the restriction from the cohomology of $K$ intersected with the invariants under $W_Q(V) = N_G(V)/V$. In our case it is clear that $E$ and $F$ are both weakly closed in $H \subset 2A_8$ (since they are the only $(\mathbb{Z}/2)^4$'s in $H$), and so we deduce that in fact

**Proposition 1.3**:
  The restriction homomorphism from $H^*(2A_8, \mathbb{F}_2)$ to the cohomology of both $E$ and $F$ maps onto the $L_3(2)$ invariants described above. ∎

  The Quillen-Venkov theorem (see [AM1]) shows that we have an exact sequence

$$0 \to Rad\ (H^*(2A_8)) \to H^*(2A_8) \xrightarrow{res_E^G \oplus res_F^G} H^*(E)^{L_3(2)} \oplus H^*(F)^{L_3(2)} \to \mathcal{M}_\mathcal{G} \to 0$$

where the term on the left is the *radical*, i.e. the ideal of nilpotent elements in the cohomology while the term on the right is the quotient ring of *double image classes* – classes which must restrict non-trivially to both $H^*(E)$ and $H^*(F)$.

In what follows we will show that $\mathcal{M}_\mathcal{G} \cong \mathbb{F}_2[d_4, d_8]$, whereas

$$Rad\ (H^*(2A_8)) \cong \mathbb{F}_2[d_4, d_8](u_3, u_7, u_9).$$

## §2. DOUBLE IMAGE CLASSES

To continue the calculation of $H^*(G)$ we now determine the term

$$\mathcal{M}_\mathcal{G} = coker\ \left(res_F^G \oplus res_E^G\right)$$

which we can identify with "double image classes" in $H^*(E) \oplus H^*(F)$ following the notions in [AM2]. We begin by observing that $H^*(E)^{L_3(2)}$ will map onto $\mathcal{M}_\mathcal{G}$. Using this projection, we can endow the vector space $\mathcal{M}_\mathcal{G}$ with a ring structure; using the other component $H^*(F)^{L_3(2)}$ we obtain the same structure, by symmetry. Next we relate the structure of $\mathcal{M}_\mathcal{G}$ to the corresponding gadget for $H = Syl_2(2A_8)$, which we denote by $\mathcal{M}_\mathcal{H}$. Using a stability argument it is easy to see that

**Proposition 2.1**:
  The inclusion $H \subset 2A_8$ induces a monomorphism of rings $\mathcal{M}_\mathcal{G} \to \mathcal{M}_\mathcal{H}$. ∎

We also need the following elementary lemma.
**Lemma 2.2**:
  Let $V = E \cap F$, then the map $H^*(E)^{L_3(2)} \oplus H^*(F)^{L_3(2)} \to H^*(V)$ factors through $\mathcal{M}_\mathcal{G}$. Furthermore the ring structure on $\mathcal{M}_\mathcal{G}$ is compatible with that in $H^*(V)$. ∎

Now we determine the structure of the double image classes as follows.

**Lemma 2.3**:
The ring $\mathcal{M}_\mathcal{G} = \mathbb{F}_2[d_4, d_8]$.

**Proof**:
The ring $\mathcal{M}_\mathcal{H}$ was computed in [AM2], and we have

$$\mathcal{M}_\mathcal{H} \cong \mathbb{F}_2[d_4, d_8](1, w_3, n_7).$$

Clearly $d_4$ and $d_8$ remain as double image classes for $G$ using 2.2.

On the other hand, since there are no elements in $\mathcal{M}_\mathcal{H}$ in dimensions 6 or 10, and the only element in $\mathcal{M}_\mathcal{H}$ having dimension 8 is $d_8$ it follows that $d_6$ is not in the double image so there are classes $d_6^1$ and $d_6^2$ in $H^*(2A_8)$ which restrict to $(d_6, 0)$, $(0, d_6)$ respectively. Likewise, there must exist classes $a_8^1$ and $a_8^2$ in $H^*(2A_8)$ which restrict to $(a_8, 0)$ and $(0, a_8)$ respectively and $a_{10}^1$, $a_{10}^2$ which restrict to $(a_{10}, 0)$ and $(0, a_{10})$ respectively. ∎

## §3. DETERMINATION OF THE RADICAL

As was observed in [AM2], $H$ contains $S = Syl_2(L_3(4))$ as an index 2 subgroup and the restriction map induces an injection of the radical $Rad\ H^*(H) \to Rad\ H^*(S)$ and hence an embedding $Rad\ H^*(2A_8) \subset Rad\ H^*(S)$. We recall that

$$Rad\ H^*(S) \cong \mathbb{F}_2[v_4, w_4](\gamma_2, \beta_2, \gamma_3, \beta_3, \alpha_5).$$

The normalizer of $S$ in $G$ is an extension of the form $S \cdot S_3$, hence $Rad\ H^*(2A_8) \subset Rad\ H^*(S)^{S_3}$. We now determine these invariants. If we write $S_3 =< T, u >$, where $T$ has order 3 and $u$ is an involution, then the explicit actions are given by (see [AM3])

$$T(v_4) = v_4, T(u_4) = u_4, T(\gamma_2) = \beta_2, T(\beta_2) = \gamma_2 + \beta_2, T(\beta_3) = \beta_3, T(\gamma_3) = \gamma_3, T(\alpha_5) = \alpha_5$$

whereas

$$u(v_4) = w_4, u(w_4) = v_4, u(\gamma_2) = \beta_2, u(\beta_2) = \gamma_2, u(\gamma_3) = \beta_3, u(\beta_3) = \gamma_3, u(\alpha_5) = \alpha_5.$$

From this action we obtain the following invariants

**Proposition 3.1**:

$$(Rad\ H^*(S))^{S_3} \cong \mathbb{F}_2[v_4 + w_4, v_4 w_4](\gamma_3 + \beta_3, \alpha_5, v_4 \gamma_3 + w_4 \beta_3).$$

∎

We can use the Lyndon–Hochschild–Serre spectral sequence for the extension

$$1 \to \mathbb{Z}/2 \to 2A_8 \to A_8 \to 1$$

to analyze the cohomology of $G$ in low dimensions. Using the results in [AMM], we in fact see that $H^i(G) = 0$ for $i = 1, 2, 5$, whereas $H^i(G) \cong \mathbb{F}_2$ for $i = 3, 4$. We also have that $H^7(G) \cong H^9(G) \cong (\mathbb{F}_2)^3$ (details are provided in the next section). On the other hand, if as in [AM2] we choose a non-nilpotent $k_5 \in H^5(H)$ which restricts to $\alpha_5$, then we may choose a class $\gamma_8 \in H^8(H)$ which restricts to $v_4 w_4$ and such that $(\gamma_8)^i k_5$ is non-nilpotent for all $i \geq 1$. The details are as follows: the products $\gamma_8^i k_5$ restrict non-trivially to the index 2 subgroup $UT_4(2)$ (upper triangular $4 \times 4$ matrices over $\mathbb{F}_2$), according to the explicit restriction data in [AM2]. However, the cohomology of this subgroup is detected on elementary abelian subgroups, from which non–nilpotence follows. Note that one can verify explicitly that $res^H_{UT_4(2)}(\gamma_4 k_5) = 0$; indeed, following the notation in [AM2] this is a consequence of the fact that $s\gamma_4 = 0$ and that $s, k_5$ restrict non–trivially to precisely the same conjugacy classes of maximal detecting elementary abelian subgroups.

Now, by looking at the image of $H^*(G)$ in the cohomology of the maximal elementary abelian subgroups, we infer that there are at most two non–nilpotent seven and nine dimensional cohomology classes respectively in $H^*(G)$, and hence there must be a seven and a nine dimensional class in the radical. Combining these fact we conclude

**Proposition 3.2**:

$$Rad\ H^*(2A_8) \cong \mathbb{F}_2[v_4 + w_4, v_4 w_4](\gamma_3 + \beta_3, v_4\gamma_3 + w_4\beta_3, \alpha_5(v_4 + w_4)).$$

∎

Summarizing, we have determined the cohomology of $2A_8$:

**Theorem 3.3**:
Let $G = 2A_8$; then we have a short exact sequence

$$0 \to \mathbb{F}_2[d_4, d_8](u_3, u_7, u_9) \to H^*(G) \stackrel{res^G_E \oplus res^G_F}{\longrightarrow} H^*(E)^{L_3(2)} \oplus H^*(F)^{L_3(2)} \to \mathbb{F}_2[d_4, d_8] \to 0$$

where

$$H^*(E)^{L_3(2)} = H^*(F)^{L_3(2)} = \mathbb{F}_2[d_4, d_6, d_7, d_8](1, a_8, a_9, a_{10}, a_{11}, a_{12}, a_{13}, a_{21}).$$

∎

§4. **THE SPECTRAL SEQUENCE OF THE CENTRAL EXTENSION**

In this section we analyze the spectral sequence associated to the central extension

$$1 \to \mathbb{Z}/2 \to G \to A_8 \to 1.$$

This allows us, initially, to determine enough differentials in the spectral sequence to show that $H^7(2A_8) \cong H^9(2A_8) \cong (\mathbb{F}_2)^3$ as claimed in the proof of 3.1. Then, using 3.3 we complete the analyis of the spectral sequence obtaining very interesting interior differentials. To begin we recall the cohomology of the alternating group $A_8$ [AMM]:

$$H^*(A_8) \cong \mathbb{F}_2[\sigma_2, \sigma_3, c_3, \sigma_4, d_6, e_6, d_7, e_7](1, x_5)/<R>$$

where the subscript denotes the degree and $R$ is the following set of relations:

$$d_6\sigma_3 = e_6\sigma_3 = \sigma_3(\sigma_2^3 + \sigma_2\sigma_4)$$
$$d_6e_7 = d_7e_6 = d_7e_7 = 0$$
$$d_7\sigma_2 = d_7\sigma_3 = d_7c_3 = d_7x_5 = 0$$
$$e_7\sigma_2 = e_7\sigma_3 = e_7c_3 = e_7x_5 = 0$$
$$x_5\sigma_3 = c_3\sigma_3 = 0$$
$$d_6e_6 = \sigma_2(\sigma_2c_3x_5 + c_3^2\sigma_4 + \sigma_2\sigma_4^2 + \sigma_2^5) + x_5c_3\sigma_4$$
$$x_5^2 + x_5\sigma_2c_3 + (d_6 + e_6)\sigma_2^2 + \sigma_4c_3^2 = 0$$

The following tables partially describe the action of the Steenrod algebra:

|        | $\sigma_2$        | $\sigma_3$        | $c_3$             |
|--------|-------------------|-------------------|-------------------|
| $Sq^1$ | $\sigma_3 + c_3$  | 0                 | 0                 |
| $Sq^2$ | $\sigma_2^2$      | $\sigma_2\sigma_3$ | $\sigma_2c_3 + x_5$ |
| $Sq^3$ | 0                 | $\sigma_3^2$      | $c_3^2$           |

|        | $\sigma_4$                                    | $x_5$                        | $d_6$                                | $d_7$       |
|--------|-----------------------------------------------|------------------------------|--------------------------------------|-------------|
| $Sq^1$ | $x_5$                                         | 0                            | $d_7 + \sigma_3(\sigma_2^2 + \sigma_4)$ | 0           |
| $Sq^2$ | $\sigma_2\sigma_4 + d_6 + e_6$                | $\sigma_2 x_5$               | —                                    | 0           |
| $Sq^3$ | $d_7 + e_7 + c_3\sigma_4 + \sigma_3\sigma_4 + \sigma_2 x_5$ | $(c_3 + \sigma_3)x_5$ | —                                    | 0           |
| $Sq^4$ | $\sigma_4^2$                                  | $\sigma_4 x_5 + c_3(d_6 + e_6)$ | —                                    | $\sigma_4 d_7$ |
| $Sq^5$ | 0                                             | $x_5^2$                      | —                                    | 0           |
| $Sq^6$ | 0                                             | 0                            | $d_6^2$                              | $d_6 d_7$   |

The action of the Steenrod squares on $e_6, e_7$ is analogous to that on $d_6, d_7$ substituting the corresponding values.

We now consider the spectral sequence of the central extension. It has

$$E_2^{*,*} \cong H^*(A_8) \otimes \mathbb{F}_2[e],$$

where $e$ represents the one dimensional generator on the fiber. We have that $d^2(e) = \sigma_2$. From the multiplicative relations above and a direct calculation using Macaulay to obtain the kernel of multiplication by $\sigma_2$ in $H^*(A_8)$, we obtain the following $E_3$-term:

$$\left\{ \begin{array}{c} \mathbb{F}_2[\sigma_3 + c_3, e_6, d_6](1, x_5)/\mathcal{R}' \oplus \mathbb{F}_2[\sigma_3 + c_3]\sigma_3 \\ \oplus \mathbb{F}_2[d_6, d_7](d_7, d_7e) \oplus \mathbb{F}_2[e_6, e_7](e_7, e_7e) \end{array} \right\} \otimes \mathbb{F}_2[\sigma_4, e^2].$$

Here, $\mathcal{R}'$ is the set of relations

$$\begin{aligned}
x_5^2 &= \sigma_4(c_3 + \sigma_3)c_3, \\
d_6 e_6 &= x_5(\sigma_3 + c_3)\sigma_4, \\
e_6 \sigma_3 &= 0, \\
d_6 \sigma_3 &= 0,
\end{aligned}$$

together with the relations for the algebra in $E_2$. Note the two nontrivial permanent cocycles $d_7 e, e_7 e$ in the $E_3^{7,1}$-term. They represent $a_8$ and $a_8'$ respectively.

Next we have $d^3(e^2) = Sq^1 \sigma_2 = c_3 + \sigma_3$. Moreover, from the form of the $E_3$-term, we see that this is the only $d^3$-differential and $E_4$ has the form

4.1
$$\left\{ \begin{array}{c} \mathbb{F}_2[e_6, d_6](1, x_5)/\mathcal{R}'' \oplus \sigma_3 \oplus \mathbb{F}_2[d_6, d_7](d_7, ed_7, e^2 d_7, e^3 d_7) \\ \oplus \mathbb{F}_2[e_6, e_7](e_7, ee_7, e^2 e_7, e^3 e_7) \end{array} \right\} \otimes \mathbb{F}_2[\sigma_4, e^4]$$

where $\mathcal{R}''$ is the pair of relations $x_5^2 = 0$, $d_6 e_6 = 0$.

There are two *possible* $d^4$-differentials:

$$\begin{aligned}
d^4(e^3 d_7) &= \epsilon d_6 x_5 \\
d^4(e^3 e_7) &= \epsilon' e_6 x_5,
\end{aligned}$$

where $\epsilon$, $\epsilon'$ are 0 or 1. We will see later that these differentials actually must occur, but, in any case there can be no $d^4$-differential in total dimension $< 10$.

Now using $Sq^2$ we see that $d^5(e^4) = x_5$. Again the relations imply that $\sigma_3 e^4$ and $x_5 e^4$ are non-trivial permanent cocycles in dimensions $\leq 9$. Thus, in total degree seven the following classes survive to $E_\infty$: $d_7, e_7, \sigma_3 e^4$, whereas in total degree nine we have $x_5 e^4, d_7 e^2, e_7 e^2$. This fills in the required information used in the computation of $H^*(G)$.

However we can now use 3.3 to determine the $E_\infty$ term in its totality, namely

**Theorem 4.2**:

In the spectral sequence for the extension $1 \to \mathbb{Z}/2 \to 2A_8 \to A_8 \to 1$ we have that there are two interior $d^4$ differentials, $d^4(e^3 d_7) = d_6 x_5$ and $d^4(e^3 d_7) = e_6 x_5$. The only remaining differential is $d^5(e^4) = x_5$, so $E_6 = E_\infty$ is isomorphic to

$$\left\{ (c_3, e^4 c_3, e^4 x_5) \oplus \left[ \begin{array}{c} \mathbb{F}_2[e_6] \oplus \mathbb{F}_2[d_6]d_6 \oplus \mathbb{F}_2[d_6, d_7](1, e, e^2, e^3 d_7)d_7 \\ \oplus \mathbb{F}_2[e_6, e_7](1, e, e^2, e^3 e_6)e_6 \end{array} \right] (1, e^4) \right\} \otimes \mathbb{F}_2[\sigma_4, e^8].$$

**Proof**:

We need to consider $d^4 : E^{7,3} \to E^{11,0}$. If there is no $d^4$ differential here, then 4.1 shows that $E_4 = E_5$ and, since $d^5(e^4) = x_5$ it follows that $E_6$ would equal $E_\infty$ and would have the form

$$\left\{ \begin{array}{c} \mathbb{F}_2[e_6, d_6](1, e^4 x_5)/\mathcal{R}'' \oplus \sigma_3(1, e^4) \oplus \mathbb{F}_2[d_6, d_7](1, e)(1, e^2)(1, e^4)d_7 \\ \mathbb{F}_2[e_6 e_7](1, e)(1, e^2)(1, e^4)e_7 \end{array} \right\} \otimes \mathbb{F}_2[\sigma_4, e^8]$$

which is clearly inconsistent with the form of the answer previously determined. The only class available to be hit by $d_7 e^3$ is $d_6 x_5$, and similarly the only class available to be hit by $e_7 e^3$ is $e_6 x_5$. It follows that these differentials must be non-trivial and $E_5$ has the form

$$\left\{ \begin{array}{c} (c_3, x_5) \oplus \mathbb{F}_2[e_6] \oplus \mathbb{F}_2[d_6]d_6 \oplus \mathbb{F}_2[d_6, d_7](1, e, e^2, e^3 d_7)d_7 \\ \mathbb{F}_2[e_6, d_7](1, e, e^2, e^3 e_7)e_7 \end{array} \right\} \otimes \mathbb{F}_2[\sigma_4, e^4].$$

The $d^5$ differential on $e^4$ now implies that $E_6$ has the form given in 4.2 which has the same Poincaré series as the answer obtained in §3. ∎

## §5. THE SYLOW SUBGROUP OF $Ly$ AND THE COHOMOLOGY OF $2S_8$

The Sylow subgroup of $Ly$ – which is also $Syl_2(2A_{10})$ – can be given as an extension in many ways. For example, as a split extension it has the form $Syl_2(J_2)\colon 2$ and $Syl_2(McL)\colon 2$, while as a non-split central extension it is given as

$$2 \triangleleft Syl_2(Ly) \xrightarrow{\pi} 2 \wr 2 \wr 2 = Syl_2(S_8).$$

A convenient description is as the extension $Syl_2(L_3(4))\colon 2^2$ where $Syl_2(L_3(4))$ is the $3\times 3$-upper triangular matrices with entries in $\mathbb{F}_4$:

$$Syl_2(L_3(4)) \;=\; UT_3(4) \;=\; \left\langle \begin{pmatrix} 1 & \alpha & \beta \\ 0 & 1 & \gamma \\ 0 & 0 & 1 \end{pmatrix} \;\middle|\; \alpha, \beta, \gamma \in \mathbb{F}_4 \right\rangle$$

and $2^2$ is generated by the two automorphisms of $Syl_2(L_3(4))$, Galois conjugation, $g$,

$$\begin{pmatrix} 1 & \alpha & \beta \\ 0 & 1 & \gamma \\ 0 & 0 & 1 \end{pmatrix}^g = \begin{pmatrix} 1 & \alpha^2 & \beta^2 \\ 0 & 1 & \gamma^2 \\ 0 & 0 & 1 \end{pmatrix}$$

and the automorphism corresponding to the automorphism of $L_3(4)$ given by $M \leftrightarrow M^{-t}$ for each $M \in GL_3(\mathbb{F}_4)$:

$$\begin{pmatrix} 1 & \alpha & \beta \\ 0 & 1 & \gamma \\ 0 & 0 & 1 \end{pmatrix}^A = \begin{pmatrix} 1 & \gamma & \beta + \alpha\gamma \\ 0 & 1 & \alpha \\ 0 & 0 & 1 \end{pmatrix}.$$

Note that $A$ is really the composite

$$m \mapsto \begin{pmatrix} 0 & 0 & 1 \\ 0 & 1 & 0 \\ 1 & 0 & 0 \end{pmatrix} m^{-t} \begin{pmatrix} 0 & 0 & 1 \\ 0 & 1 & 0 \\ 1 & 0 & 0 \end{pmatrix}$$

and thus extends to an automorphism of all of $GL_3(\mathbb{F}_4)$ leaving the subgroup of upper triangular matrices invariant. Note also that $UT_3(4)$ contains exactly two copies of $2^4$: the subgroups

$$2_I^4 \;=\; \left\langle \begin{pmatrix} 1 & \alpha & \beta \\ 0 & 1 & 0 \\ 0 & 0 & 1 \end{pmatrix} \right\rangle, \qquad 2_{II}^4 \;=\; \left\langle \begin{pmatrix} 1 & 0 & \alpha \\ 0 & 1 & \beta \\ 0 & 0 & 1 \end{pmatrix} \right\rangle.$$

We also denote the generators of the center of $UT_3(4)$ as

$$T = \begin{pmatrix} 1 & 0 & \zeta_3 \\ 0 & 1 & 0 \\ 0 & 0 & 1 \end{pmatrix}, \quad Z = \begin{pmatrix} 1 & 0 & 1 \\ 0 & 1 & 0 \\ 0 & 0 & 1 \end{pmatrix}.$$

Note that the center of $Syl_2(Ly) = UT_3(4){:}\langle g, A\rangle$ is the single $\mathbb{Z}/2 = \langle Z\rangle$. The subgroup $UT_3(4){:}\langle g\rangle = Syl_2(McL)$ and is also $Syl_2(2A_8)$. It has center $\langle Z\rangle$ and exactly two subgroups isomorphic to $2^4$, $2^4_I$, and $2^4_{II}$, both of which are normal.

The subgroup $UT_3(4){:}\langle gA\rangle = Syl_2(J_2) = Syl_2(J_3)$ and here the two $2^4$'s are fused. Finally, in the full group, $Syl_2(Ly) = Syl_2(2S_8) = Syl_2(2A_{10})$, there are again only the two $2^4$ subgroups above, but they are conjugate. However, there are a number of new $2^3$'s, for example $\langle g, A, Z\rangle$ and $\langle A, T, Z\rangle$. In all three of the groups $2S_8 \subset 2A_{10} \subset Ly$ there are precisely two conjugacy classes of maximal elementary abelian subgroups, and we can take as representatives the $2^4$ above and the $2^3$ defined by $\langle A, T, Z\rangle$. This of course immediately implies that both of these elementary abelian groups form part of weakly closed systems in $2S_8 \subset 2A_{10}$ and $2S_8 \subset Ly$.

Now we consider the extension $2A_8{:}2 = 2S_8 \subset 2A_{10}$. The action of the extending $\mathbb{Z}/2$ exchanges the two $2^4$'s, consequently, on passing to cohomology it exchanges the generators corresponding to $(d_6, 0)$, $(0, d_6)$, $(d_7, 0)$, $(0, d_7)$ and the similar elements corresponding to the terms $(a_i, 0)$, $(0, a_i)$, while it and fixes the remaining generators. It follows that the $E_2$-term of the resulting Lyndon-Hochschild-Serre spectral sequence is

$$\mathbb{F}_2[d_4, d_6, d_7, d_8](1, a_8, a_9, a_{10}, a_{11}, a_{12}, a_{13}, a_{21}) \bigoplus \mathbb{F}_2[d_4, d_8, w](w, u_3, u_7, u_9)$$

where $w$ represents the cohomology class dual to the extending $\mathbb{Z}/2$.

**Lemma 5.1:**
There are two conjugacy classes of maximal elementary 2-groups in $2S_8$, $2A_{10}$, and at most two in $Ly$, the $2^4$ discussed above which is weakly closed in $Syl_2(2S_8) \subset 2S_8$, and one of the new $2^3$'s – the choice we will take in the following is $\langle A, T, Z\rangle$ – which is weakly closed in $2S_8 \subset 2A_{10}$. The Weyl groups of $\langle A, T, Z\rangle$ in $2S_8$ and $2A_{10}$ are $D_8$ and $S_4$ respectively.

**Proof**: We start by verifying the following claim.

**Claim:** The maximal elementary abelian 2-groups in $2A_{10}$, $2S_8$ are in one to one correspondence (via the quotient maps) with the elementary abelian 2-groups in $A_{10}$, $S_8$ whose nonidentity elements are all products of four disjoint transpositions and are maximal with respect to this condition.

To see this, one should note that $2A_5$ has no non-central involutions, and both $A_8$ and $A_{10}$ have two conjugacy classes of involutions, one which lifts to an involution while the other lifts to an element of order four [Co]. We also have the commuting diagram

$$\begin{array}{ccccc} 2A_5 & \hookrightarrow & 2A_8 & \hookrightarrow & 2A_{10} \\ \downarrow & & \downarrow & & \downarrow \\ A_5 & \hookrightarrow & A_8 & \hookrightarrow & A_{10} \end{array}$$

which tells us that we must have the standard embedding of $A_5$ in $A_{10}$. Since the non-central involutions of $A_{10}$ must project onto involutions not contained in $A_5$, they must have cycle type $(1,2)(3,4)(5,6)(7,8)$.

If we consider $S_8 \subset A_{10}$ generated by $(1,2)(9,10)$, $(2,3)(9,10)$, ... $(7,8)(9,10)$, then we see that the remarks above apply to $S_8$ as well.

**Claim:** The maximal elementary abelian 2-groups in $S_n$ are (up to conjugacy) of the form $V_1^{i_1} \times \cdots \times V_r^{i_r}$, where $V_j \cong 2^j$ and the permutation action is the regular representation. Here, the number of letters permuted by the subgroup is $2i_1 + \cdots + 2^r i_r$, and this number is either $n$ or $n-1$ (see e.g., [AM1], Chapter VI).

Observe that the subgroups $V_i$ described above are contained in the alternating group if $i \geq 2$. Furthermore, note that the subgroup $E_r$ of even permutations contained in $V_1^r$ is not maximal in $A_{2r}$ if $r = 2$, since the even permutations in $V_1 \times V_1$ are strictly contained in $V_2$. On the other hand if $r \geq 3$ then $E_r$ is a maximal elementary abelian 2-group in $A_{2r}$, and is the unique such group with $r$ 2-element orbits.

Altering the arguments for the symmetric group to reflect these remarks we obtain the following:

**Proposition 5.2:**

*Every maximal elementary abelian 2-group in $A_n$ is conjugate to one of the form $E_{i_1} \times V_2^{i_2} \times \cdots \times V_r^{i_r}$, where $V_j \approx 2^j$ and the permutation action is the regular representation of $2^j$. Here, the number of letters permuted by the subgroup is $2i_1 + 4i_2 + \cdots + 2^r i_r$, and this number is $n$, $n-1$, $n-2$, or $n-3$.*

To complete the classification of the conjugacy classes of elementary abelian 2-groups in $A_n$, we need the following remark.

REMARK: Suppose that $8n = 2^3 i_1 + 2^4 i_4 + \cdots + 2^r i_r$. Then a conjugacy class of maximal elementary abelian 2-groups in $S_{8n}$ or $S_{8n+1}$ of the form $V_3^{i_3} \times V_4^{i_4} \times \cdots \times V_r^{i_r}$ (which is of course also maximal in $A_{8n}$, $A_{8n+1}$, $A_{8n+2}$, or $A_{8n+3}$) splits into two conjugacy classes of elementary abelian 2-groups in $A_{8n}$ or $A_{8n+1}$. The other types of maximal elementary abelian 2-groups in $A_m$ are conjugate if and only if they have the same form.

The essential point in the proof of this remark is the fact that the number of conjugates of a subgroup $E \subset G$ is the index $[G : N_G(E)]$ of the normalizer. Applying this to the case $E \subset A_m$, we see that the number of conjugates of a $E$ in $S_m$ is either the same as the number of conjugates of $E$ in $A_m$, or twice as many, depending on whether there is a transposition of $S_m$ normalizing $E$. Having made this point, our remark reduces to the fact that if $E = V_3^{i_3} \times V_4^{i_4} \times \cdots \times V_r^{i_r}$, there is no transposition of $S_{8n}$ (or $S_{8n+1}$) normalizing $E$, and that there is one in the other cases.

From what we have shown it follows that (up to conjugacy) the maximal elementary abelian 2-groups in $S_8$ are $V_3$, $V_2^2$, $V_2 \times V_1^2$, and $V_1^4$. The conjugacy classes of maximal elementary abelian subgroups of $A_{10}$ are represented by $V_3$, $V_2^2$, $V_2 \times E_3$, and $E_5$.

Looking at the subgroups of the groups from the lists above, the only ones whose elements satisfy the necessary condition on cycle type are the conjugacy classes of

$$\tilde{V}_3 = <(1,2)(3,4)(5,6)(7,8), (1,3)(2,4)(5,7)(6,8), (1,5)(2,6)(3,7)(4,8)>$$

and

$$\tilde{M}_3 = <(1,3)(2,4)(5,6)(9,10), (1,4)(2,3)(7,8)(9,10)>.$$

Then $M_3$ and $V_3$, the lifts of $\tilde{M}_3$ and $\tilde{V}_3$, represent the maximal elementary abelian 2-groups of both $2A_{10}$ and $2S_8$.

Finally, we need to check the Weyl groups of $M_3$ in $2S_8$ and $2A_{10}$ are as stated. Clearly, $V_2 = \langle (1,3)(2,4), (1,4)(2,3) \rangle$ and $(1,2)(5,7)(6,8)$ normalize this group, and these together generate a $D_8$. But since $(9,10)$ must be fixed by the elements of $2S_8$ it follows that there can be no element of order three which normalizes it in $2S_8$. But the element of order three in $2A_{10}$ which cyclically permutes the transpositions $(5,6)$, $(7,8)$, and $(9,10)$ while simultaneously taking $(1,3)(2,4)$ to $(1,2)(3,4)$ also normalizes $M_3$ and the proof of 5.1 is complete. ∎

We have the following calculation of the invariant subring for $\langle A, T, Z \rangle$ under the action of the $D_8$-Weyl group for $2S_8$:

**Lemma 5.3**:
$\mathbb{F}_2[w,t,z]^{D_8} = \mathbb{F}_2[w, t(t+w), d_4]$ where $w$ is dual to the extending element $A = (1,2)(5,6)(7,8)$, $t$ is dual to $(1,2)(3,4)(5,8)(6,7)$ corresponding to $T$, and $z$ is dual to the central element, $Z$. Also,
$$d_4 = z^4 + z^2 d_2 + z d_3 + d_2^2$$
is the generating element of the Dickson algebra, with $d_2 = w^2 + tw + t^2$, $d_3 = wt(w+t)$.

**Proof**:
We start with the action of the element of order two which has the form $z \mapsto z+t$, $t \mapsto t$, $w \mapsto w$. This give the invariant subring
$$\mathbb{F}_2[w, t, z(z+t)].$$

Now consider the action of the second element of order two fixing $t$, $w$, $z \mapsto z+w$. Then $z(z+t) \mapsto (z+w)(z+t+w) = z(z+t) + w(t+w)$ and the invariant subring is clearly $\mathbb{F}_2[w, t, z(z+t)(z+w)(z+t+w)]$. Expanding out this last term is $z^4 + z^2 d_2 + z d_3$. Finally, the last element fixes $w$, $z$ and takes $t$ to $t+w$. Applying this, we clearly get the asserted ring of invariants. ∎

Now here is the main result of this section.

**Theorem 5.4**:
In the spectral sequence above for $H^*(2S_8)$ we have that $E_2 = E_\infty$. Moreover, the restriction map $H^*(2S_8) \to H^*(2^4)^{L_3(2)} \oplus H^*(2^3)^{D_8}$ is injective in cohomology. Moreover, the quotient by the image has the form $\mathbb{F}_2[d_4, d_8]$ with $d_4$ being the image of $d_4$ from $H^*(2^4)$ and $(w^2 + t(w+t))^2$ from $H^*(2^3)$ while $d_8$ is the image of $d_8$ and $d_4^2$, respectively.

**Proof**:
We have ([AM1] p. 212) that $H^*(S_8) = \mathbb{F}_2[\sigma_1, \sigma_2, \sigma_3, \sigma_4, c_3, d_6, d_7](x_5)/R$ where $R$ is the set of relations
$$d_6 \sigma_1 = d_6 \sigma_3 = 0$$
$$d_7 \sigma_1 = d_7 \sigma_2 = d_7 \sigma_3 = d_7 c_3 = d_7 x_5 = 0$$
$$x_5 \sigma_3 + c_3 \sigma_4 \sigma_1 = 0$$
$$c_3(\sigma_3 + \sigma_1 \sigma_2) + \sigma_1 x_5 = 0$$
$$x_5^2 + x_5 \sigma_2 c_3 + d_6 \sigma_2^2 + \sigma_4 c_3^2 = 0.$$

We also have the commutative diagram

$$\begin{array}{ccccc}
\mathbb{Z}/2 & \xrightarrow{=} & \mathbb{Z}/2 & \xrightarrow{=} & \mathbb{Z}/2 \\
\downarrow & & \downarrow & & \downarrow \\
2^3 & \hookrightarrow & 2^3 \cdot 2^2 & \hookrightarrow & 2S_8 \\
\downarrow & & \downarrow & & \downarrow \\
W = 2^2 & \xrightarrow{\Delta} & V_1^2 \times K & \hookrightarrow & S_8,
\end{array}$$

where $\Delta$ is the diagonal map, (identifying both $V_1^2$ and $K$ with $2^2$ and noting that the generators of $W$ go diagonally to generators of the first and second $2^2$'s). From this it is direct to determine the image of the restriction map from the generators of $H^*(S_8)$ to $H^*(2^3)$ as follows

$$\begin{cases}
\sigma_1 \mapsto \sigma_1 \otimes 1 \mapsto w \\
\sigma_3 \mapsto \sigma_1 \otimes d_2 \mapsto w^3 + d_3 \\
c_3 \mapsto 1 \otimes d_3 \mapsto d_3 \\
\sigma_4 \mapsto \sigma_2 \otimes d_2 \mapsto t(t+w)d_2 \\
x_5 \mapsto \sigma_2 \otimes d_3 \mapsto t(t+w)d_3 \\
d_6 \mapsto 0 \\
d_7 \mapsto 0.
\end{cases}$$

This shows that the subalgebra $\mathbb{F}_2[w, d_2^2](1, d_3)$ in $H^*(2^3)$ is contained in the image from the cohomology of $H^*(2S_8)$ as the subring coming from the composition

$$H^*(S_8) \to H^*(2S_8) \xrightarrow{res} H^*(2^3).$$

Also, by the argument of [AM3], §6, there is a real representation of $Ly$ so that the Stiefel-Whitney class $w_8$ restricts to $d_8$ in $H^*(2^4)$ and $d_4^2 \in H^*(\langle Z, W\rangle)$. It follows that the image of $res\colon H^*(2S_8)\to H^*(\langle ZW\rangle)$ contains the subalgebra $\mathbb{F}_2[w, d_2^2, d_4^2](1, d_3)$.

We now use the Lyndon-Hochschild-Serre spectral sequences of the central extensions $2 \triangleleft 2^3 \xrightarrow{\pi} W$ and $2 \triangleleft 2S_8 \xrightarrow{\pi} S_8$ together with naturality to show that the remaining classes $a_7$ and $X_9$ also restrict non-trivially to $H^*(\langle Z, W\rangle)$.

In the spectral sequence for $2S_8$ we write $E_2 = H^*(S_8) \otimes \mathbb{F}_2[z]$, while in the spectral sequence for $2^3$ we have

$$\begin{aligned}
E_2 &= H^*(W) \otimes \mathbb{F}_2[z] \\
&= \mathbb{F}_2[w, t] \otimes \mathbb{F}_2[z].
\end{aligned}$$

Moreover, the restriction map $res\colon E_2(2S_8)\to E_2(2^3)$ has the form $res \otimes id$, and is natural with respect to differentials. Of course, the differentials in the spectral sequence for $2^3$ are all trivial, which implies that no element in the spectral sequence for $H^*(2S_8)$ which restricts non-trivially to the sequence for $2^3$ *can be in the image of a differential* in the spectral sequence for $2S_8$.

The differential $d^2$ is given by $d^2(z) = \sigma_2 + \sigma_1^2$ while the basic higher differentials are given by

$$\begin{aligned}
d^3(z^2) &= \sigma_1\sigma_2 + \sigma_3 + c_3 \\
d^5(z^4) &= x_5 + c_3\sigma_1^2 + \sigma_1\sigma_4.
\end{aligned}$$

An easy calculation shows that $d^5(x_5 z^4) \in H^{10}(S_8)$ is in the image of the differentials $d^2$ and $d^3$, and a similar calculation shows the same thing for $d^5(c_3 z^4)$. Consequently, since these classes restrict non-trivially to the spectral sequence for $H^*(2^3)$, it follows that they survive to $E_\infty$ and represent non-trivial classes in $H^*(2S_8)$.

But the argument above actually shows more. It shows that, up to filtration, the class represented by $x_5 z^4$ restricts to $res(x_5) z^4$, i.e., the image of restriction on this class has the form
$$t(t+w) d_3 z^4 + v_6(t,w) z^3 + v_7(t,w) z^2 + v_8(t,w) z + v_9(t,w).$$
On the other hand the image of this class lies in the invariant subring, so it must have the form $t(t+w) d_3 d_4 + v_9(t,w)$ with $v_9(t,w)$ also in the invariant subring. However, it is direct to check that every element in dimension nine in $\mathbb{F}_2[w, t(t+w)]$ is already in the image of restriction from $H^*(S_8)$, and $t(t+w) d_3 d_4$ is in the image of restriction from $H^*(2S_8)$.

A similar argument with $c_3 z^4$ shows that $d_3 d_4$ is in the image of restriction from $H^*(2S_8)$. But the smallest sub-ring of $H^*(2^3)$ containing these elements and $res(H^*(S_8))$ has the form
$$\mathbb{F}_2[w, d_2^2, d_4^2](1, d_3, d_3 d_4, t(t+w) d_3 d_4).$$
This completes the proof. ∎

**Corollary 5.5**:

The group $\langle A, T, Z \rangle \cong M_3$ is weakly closed in $2S_8 \subset Ly$. Also, the Weyl groups of $\langle A, T, Z \rangle$ in $Ly$ is $L_3(2)$.

**Proof**: Note that from the work on $2S_8$ we already know that $H^*(Ly)$ is detected by restriction to $H^*(2^4)^{A_7} \oplus H^*(2^3)$. If the $2^3$ fused into $2^4$ then $H^*(Ly)$ would be detected by $H^*(2^4)^{A_7}$ alone. But from [AM3], §6, we have

**Proposition 5.6**:

There is a real representation of $Ly$ so that the Stiefel-Whitney class $w_8$ restricts to $d_8$ in $H^*(2^4)$ and to $d_4^2 \in H^*(\langle A, T, Z \rangle)$.

Since there are no other invariants in dimensions 7 or 9 it follows that $d_8$ must be the restriction of a torsion free cohomology class and this is impossible. Consequently the $2^3$ cannot fuse into the $2^4$.

Similarly, we show that the Weyl group of this $2^3$ is $L_3(2)$. Indeed, since the part of the proof of 5.1 already finished shows that the Weyl group is at least $\mathcal{S}_4$ which is maximal in $L_3(2)$, it follows that the Weyl group is either $\mathcal{S}_4$ or $L_3(2)$. But if it were $\mathcal{S}_4$ then weak closure would imply that the restriction to $H^*(2^3)$ is the same as the restriction of $H^*(2A_{10})$. On the other hand both $d_3$ and $d_2^2$ are in this image, but $d_2^2$ is a double image class, and hence $d_4$ would have to be in the invariant subring $H^*(2^4)^{A_7}$. This is impossible, so it follows that the Weyl group is $L_3(2)$. ∎

Summarizing we have the following result which reduces the calculation of $H^*(2A_{10})$ and $H^*(Ly)$ to direct calculations with rings of invariants.

**Corollary 5.7**:

$H^*(2A_{10})$ and $H^*(Ly)$ are both detected by restriction to the direct sum $H^*(2^4) \oplus H^*(\langle A, T, Z \rangle)$. In particular, since
$$\langle Z, A, T \rangle = \langle Z, W \rangle \subset 2S_8 \subset \begin{cases} 2A_{10} \\ Ly \end{cases}$$

are both weakly closed systems it follows that the image of restriction from $H^*(2A_{10})$ to $H^*(\langle Z,W\rangle)$ is given as the intersection

$$\mathbb{F}_2[w, d_2^2, d_4^2](1, d_3, d_3d_4, t(t+w)d_3d_4) \cap H^*(\langle Z,W\rangle)^{\mathcal{S}_4}$$

while the image of restriction from $H^*(Ly)$ is given as

$$\mathbb{F}_2[w, d_2^2, d_4^2](1, d_3, d_3d_4, t(t+w)d_3d_4) \cap \mathbb{F}_2[d_4, d_6, d_7].$$

## §6. THE COHOMOLOGY OF $2A_{10}$

We will now make explicit the cohomology of $2A_{10}$ using the results in §5. First we recall that the $\mathcal{S}_4$–invariants can be computed from the usual symmetric invariants and in our notation can be expressed as

$$\mathbb{F}_2[w, t, z]^{\mathcal{S}_4} = \mathbb{F}_2[d_2, d_3, d_4].$$

Hence we immediately obtain the classes $d_3$, $d_2^2$, $d_3d_4$, and $d_4^2$ in the intersection. Moreover note the relations

$$w^2 d_3 + w d_2^2 + w^5 = d_2 d_3, \qquad w^2 + t(t+w)d_3d_4 = d_2d_3d_4.$$

From this we see that the classes $d_2d_3$ and $d_2d_3d_4$ are also in the intersection, and hence it is at least as big as

$$\mathbb{F}_2[d_3, d_2^2, d_4^2](1, d_2d_3, d_3d_4, d_2d_3d_4).$$

This is in fact the totality of the intersection; complete details for this computation are provided in an appendix. We can describe the cohomology of $2A_{10}$ as follows:

**Theorem 6.1**:
    The cohomology of $2A_{10}$ is described by a short exact detection sequence of the form

$$0 \to H^*(2A_{10}) \to H^*((\mathbb{Z}/2)^4)^{L_3(2)} \oplus \mathbb{F}_2[d_3, d_2^2, d_4^2](1, d_2d_3, d_3d_4, d_2d_3d_4) \to \mathbb{F}_2[h_4, h_8] \to 0$$

where $h_4, h_8$ correspond to the double image classes $(d_4, d_2^2)$ and $(d_8, d_4^2)$. ∎

## §7. THE COHOMOLOGY OF THE LYONS GROUP

Now we turn to the determination of the cohomology of the sporadic group $Ly$. From our previous results it will be completely determined once we have

$$\mathbb{F}_2[d_3, d_2^2, d_4^2](1, d_2d_3, d_3d_4, d_2d_3d_4) \cap \mathbb{F}_2[d_4, d_6, d_7].$$

**Lemma 7.1**:
    The intersection above is given as

$$\mathbb{F}_2[d_4^2, d_6^2, d_7](1, d_4d_7, d_6d_7, d_4d_6d_7) = \mathbb{F}_2[d_4^2, d_6^2] \oplus \mathbb{F}_2[d_4, d_6, d_7]d_7.$$

**Proof**:

We can write the ring

$$\mathbb{F}_2[d_3, d_2^2, d_4^2](1, d_2d_3, d_3d_4, d_2d_3d_4) = \mathbb{F}_2[d_3, d_2^2, d_4^2](1, d_2d_3) \oplus \mathbb{F}_2[d_2, d_3, d_4^2]d_3d_4.$$

On the other hand we have the decompositions

$$\begin{aligned}
d_6 &= d_2d_4 + d_2^3 + d_3^2, \\
d_7 &= d_3d_4 + d_2^2d_3,
\end{aligned}$$

so we can modifiy the decomposition above to the new decomposition

$$\mathbb{F}_2[d_3, d_2^2, d_4^2](1, d_2d_3) \oplus \mathbb{F}_2[d_2, d_3, d_4^2]d_7.$$

But we can further write the first summand in the form

$$\mathbb{F}_2[d_2^2, d_3^2, d_4^2](1, d_3, d_2d_3, d_2d_3^2).$$

Using this we can prove the

**Sublemma 7.2**:

*The intersection of the ring above with $\mathbb{F}_2[d_4, d_6]$ is $\mathbb{F}_2[d_4^2, d_6^2]$.*

**Proof**:

Clearly, the intersection is $\mathbb{F}_2[d_4, d_6] \cap \mathbb{F}_2[d_2^2, d_3^2, d_4^2](1, d_3, d_2d_3, d_2d_3^2)$. But to get a term involving $d_6^i$ we must have a term involving $d_2^i d_4^i$. Moreover, it cannot have any $d_3$'s involved. Hence the intersection must lie in $\mathbb{F}_2[d_2^2, d_3^2, d_4^2]$. But this intersection with $\mathbb{F}_2[d_4, d_6, d_7]$ is $\mathbb{F}_2[d_4^2, d_6^2, d_7^2]$. ∎

Note that we can replace $\mathbb{F}_2[d_2^2, d_3, d_4^2](d_3d_4^2, d_2d_3d_4^2)$ by $\mathbb{F}_2[d_2^2, d_3^2, d_4^2]d_7d_4$, and this summed with $\mathbb{F}_2[d_2, d_3, d_4^2]d_7$ is $\mathbb{F}_2[d_2, d_3, d_4]d_7$. It follows that the intersection contains $\mathbb{F}_2[d_4, d_6, d_7]d_7$ and since the quotient of $\mathbb{F}_2[d_4, d_6, d_7]$ by this ideal is $\mathbb{F}_2[d_4, d_6]$ it follows that the intersection is given as

$$\mathbb{F}_2[d_4, d_6] \cap \mathbb{F}_2[d_3, d_2^2, d_4^2](1, d_2d_3, d_3d_4, d_2d_3d_4)$$

direct summed with the ideal $(d_7)$. ∎

Now we are able to describe the cohomology of $Ly$. Recall that the Weyl group of the $2^4$ in $Ly$ is $\mathcal{A}_7$, and that the ring of $\mathcal{A}_7$-invariants has been determined in [AM4] as

$$H^*(2^4)^{\mathcal{A}_7} = \mathbb{F}_2[d_8, d_{12}, d_{14}, d_{15}](1, x_{18}, x_{20}, x_{21}, x_{24}, x_{25}, x_{27}, x_{45}).$$

**Theorem 7.3**:

*The cohomology of $Ly$ is given by the short exact detection sequence*

$$0 \longrightarrow H^*(Ly) \longrightarrow H^*(2^4)^{\mathcal{A}_7} \oplus \mathbb{F}_2[d_4^2, d_6^2, d_7](1, d_4d_7, d_6d_7, d_4d_6d_7) \longrightarrow \mathbb{F}_2[h_8, h_{12}] \longrightarrow 0$$

*where $h_8$ and $h_{12}$ correspond to the double image classes $(d_8, d_4^2)$ and $(d_{12}, d_6^2)$.* ∎

**Remark**: We would like to point out that from the spectral sequences for computing $H^*(2A_{10})$ and $H^*(2S_8)$ it is apparent that their cohomology rings have depth equal to three i.e. there is a regular sequence of length three. On the other hand, an explicit computer algebra verification shows that the centralizers of the $(\mathbb{Z}/2)^3$ subgroups in these groups are either the $(\mathbb{Z}/2)^3$ themselves or a $(\mathbb{Z}/2)^4$. A theorem due to Carlson [C] indicates that a depth three cohomology ring must be detected on the centralizers of rank three elementary abelian subgroups. Of course, currently, the only way to verify the conditions of Carlson's result is to complete the calculations of the cohomology groups since differentials in the spectral sequences can lower ranks. But it would be very interesting if one could find independent methods for showing that these conditions are satisfied.

## APPENDIX

Let $k$ be the field of two elements and $k[d_4, d_6, d_7]$ the Dickson algebra of $GL_3(k)$-invariants in the polynomial ring $k[w, t, z]$. Let $d_2, d_3$ be the generators of the $GL_2(k)$-invariants in $k[w, t]$. In this appendix we will show that

$$k[w, d_2^2, d_4^2](1, d_3, d_3 d_4, t(t+w)d_3 d_4) \cap k[d_2, d_3, d_4] = k[d_2^2, d_3, d_4^2](1, d_2 d_3, d_3 d_4, d_2 d_3 d_4).$$

To begin we observe that since $d_2 = w^2 + t(t+w)$, it must follow that

$$k[w, d_2^2, d_4^2](1, d_3, d_3 d_4, t(t+w)d_3 d_4) = k[w, d_2^2, d_4^2](1, d_3, d_3 d_4, d_2 d_3 d_4).$$

We will require the following

**Lemma A1**

The subalgebra $k[w, d_2, d_3]$ of $k[w, t]$ is a free $k[d_2, d_3]$-module with basis $\{1, w, w^2\}$.
**Proof:** This is a special case of a more general fact about polynomial invariants (see [Mit]). However, for the sake of completeness we indicate a proof here: first, from the usual lemma on integral extensions, we have that $\{1, w, w^2\}$ generates $k[w, d_2, d_3]$ as a $k[d_2, d_3]$-module. So it is enough to prove that there is no relation $p + qw + rw^2 = 0$ with $p, q, r \in k[d_2, d_3]$. Such a relation would be an equation in the polynomial ring $k[w, t]$ and therefore remains true after applying an automorphism of that ring. So, choosing $\varphi, \theta \in GL_2(k)$ such that $\varphi(w) = t, \theta(w) = w + t$, we have the relations $p + qt + wt^2 = 0$ and $p + q(w + t) + r(w^2 + t^2) = 0$. Adding gives $p = 0$. Since $k[w, t]$ is a domain, it follows from $qw + rw^2 = 0$ that $q + rw = 0$. The same argument gives $q = 0$ and it then follows that $r = 0$. ∎

It follows from the lemma that $k[d_2, d_3, d_4](1, w, w^2)$ is a free $k[d_2, d_3, d_4]$-module with basis $\{1, w, w^2\}$. Since $k[d_2, d_3, d_4]$ is a free module over $k[d_2^2, d_3, d_4^2]$ with basis $\{1, d_2, d_4, d_2 d_4\}$, it follows that: $k[d_2, d_3, d_4](1, w, w^2)$ is a free module over $k[d_2^2, d_3, d_4^2]$ with basis the twelve elements obtained from the product $\{1, w, w^2\} \cdot \{1, d_2, d_4, d_2 d_4\}$, i.e.

$$1, d_2, d_4, d_2 d_4, w, w d_2, w d_4, w d_2 d_4, w^2, w^2 d_2, w^2 d_4, w^2 d_2 d_4.$$

We will call the free $R$-module mentioned above $M$. We use the fact that if $A \subset B \subset C$ are rings and $C$ is a finitely-generated free $B$-module with basis $\{c_i\}$ and $B$ is a finitely-generated free $A$-module with basis $\{b_j\}$, then $C$ is a finitely-generated free $A$-module with basis $\{b_j c_i\}$.

In what follows we will need the following formal lemma on the same theme.

**Lemma A2:**

If $\alpha, \beta$ are integral elements over $R$ and $\{\alpha_1, \ldots, \alpha_n\}$ generate $R[\alpha]$ as an $R$-module, while $\{\beta_1, \ldots, \beta_m\}$ generate $R[\beta]$ as an $R$-module, then $\{\alpha_i \beta_j\}$ is a set of $R$-module generators for $R[\alpha, \beta]$. ∎

The proof follows from regarding $R[\alpha, \beta]$ as $R[\alpha][\beta]$, which makes sense since $\beta$ integral over $R$ implies $\beta$ integral over $R[\alpha]$. We note that the system of generators in this lemma may be very redundant.

**The Intersection**

Let $R = k[d_2^2, d_3, d_4^2]$. It follows from the remarks in the previous section that $w$ is integral over $R$. In fact, $X^6 + d_2^2 X^2 + d_3^2 = 0$ is an integral equation for $w$ over $R$. It follows that $k[w, d_2^2, d_4^2](1, d_3)$ is an integral extension of $R$ with $\{1, w, \ldots, w^5\}$ as a set of $R$-module generators. Furthermore $d_3 d_4$, $d_2 d_3 d_4$ are integral over $R$, with equations $X^2 + d_3^2 d_4^2$ and $X^2 + d_2^2 d_3^2 d_4^2$ repectively.

It follows from Lemma A2 that the 24 elements given by all triple products from the set $\{1, w, \ldots, w^5\} \cdot \{1, d_3 d_4\} \cdot \{1, d_2 d_3 d_4\}$ generate

$$k[w, d_2^2, d_4^2](1, d_3, d_3 d_4, d_2 d_3 d_4)$$

as an $R$-module. We denote $k[w, d_2^2, d_4^2](1, d_3, d_3 d_4, d_2 d_3 d_4)$ by $U$ when regarding it as an $R$-module. When $k[d_2, d_3, d_4]$ is regarded as an $R$-module, it will be denoted by $V$. Notice that $U$ and $V$ are $R$-submodules of $M$. In the notation just described, the purpose of this note is to find generators for $U \cap V$.

To begin describing a simpler generating set for the $R$-module $U$, let us first examine the subset $\{1, w, \ldots, w^5\}$ of the generators of $U$. Since $w^3 + d_2 w + d_3 = 0$ and $d_3 \in R$, the $R$-submodule of $U$ generated by $\{1, w, w^2, w^3\}$ is the same as that generated by $\{1, w, w^2, d_2 w\}$. Similar reasoning, plus the equation $w^5 + d_2^2 w + d_2 d_3 = 0$, gives that the $R$-submodule generated by $\{1, w, \ldots, w^5\}$ is the same as that generated by $\{1, w, w^2, w d_2, w^2 d_2, d_2 d_3\}$.

Thus we can replace the first set in the 24 generators remark above with the set $\{1, w, w^2, w d_2, w^2 d_2, d_2 d_3\}$.

As $U$ is a submodule of the free module $M$ defined in the previous section we can write our new generating set in terms of the basis previously described from the free module $M$. To make it clear when we are regarding a polynomial as one of the free generators of $M$, we will represent the free generators in boldface. Thus $M$ is the free $R$-module on generators

$$\mathbf{1}, \mathbf{d_2}, \mathbf{d_4}, \mathbf{d_2 d_4}, \mathbf{w}, \mathbf{w d_2}, \mathbf{w d_4}, \mathbf{w d_2 d_4}, \mathbf{w^2}, \mathbf{w^2 d_2}, \mathbf{w^2 d_4}, \mathbf{w^2 d_2 d_4}.$$

Multiplying out the 24 triple products yields the following generators of $U$ (regarded here as elements of $M$):

$$1, w, w^2, wd_2, w^2d_2, d_3d_2, d_3d_4, d_3wd_4, d_3w^2d_4, d_3wd_2d_4, d_3w^2d_2d_4,$$

$$d_3d_2d_4, d_3d_2d_4, d_3wd_2d_4, d_3w^2d_2d_4, d_2^2d_3wd_4, d_2^2d_3w^2d_4, d_2^2d_3^2d_4,$$

$$d_3^2d_4^2d_2, d_3^2d_4^2d_2w, d_3^2d_4^2d_2w^2, d_3^2d_4^2d_2^2w, d_3^2d_4^2d_2^2w^2, d_2^2d_3^3d_4^2 1$$

On the other hand the $R$-submodule $V$ of $M$ is much simpler, being generated by $\{1, d_2, d_4, d_2d_4\}$. Since every generator of $U$ in the long list above is an $R$-multiple of a generator of $M$, we need only consider what $R$-multiples of the four generators $1, d_2, d_4, d_2d_4$ occur.

More explicitly, we are quoting the following:

**Lemma A3**

Let $R$ be a ring and $M$ a free $R$-module with basis $\{e_1, \ldots, e_n\}$. Let $V$ be the submodule generated by $\{e_1, \ldots, e_k\}$ for some $k < n$. Let $U$ be an $R$-module generated by $\{r_{ij} \cdot e_i\}$. Then

$$U \cap V = J_1 \oplus \cdots \oplus J_k \oplus 0 \oplus \cdots \oplus 0,$$

where $J_t$ is the ideal of $R$ generated by $\{r_{tj}\}$ (for fixed $t$, but $j$ varying). ∎

In any case, a brief examination of the list shows that the relevant multiples, or ideals as in the lemma are $(1), (d_3), (d_3), (d_3)$. Thus, the $R$-submodule $U \cap V$ of $M$ is generated by $1, d_3d_2, d_3d_4, d_3d_2d_4$. Now all the $R$-modules involved are actually rings, and a set of $R$-module generators is also a set of ring generators, so we have shown what we claimed at the outset, namely that:

$$k[w, d_2^2, d_4^2](1, d_3, d_3d_4, t(t+w)d_3d_4) \cap k[d_2, d_3, d_4] = k[d_2^2, d_3, d_4^2](1, d_2d_3, d_3d_4, d_2d_3d_4).$$